\documentclass[11pt]{article}
\usepackage[fleqn]{amsmath}
\usepackage{alltt}
\usepackage{amssymb}
\newtheorem{thm}{Theorem}

\def\x{$\hfill\rlap{$\sqcup$}\sqcap$\bigskip}

\begin{document}

\title{Counting the Number of Site Swap Juggling Patterns with Respect to Particular Ceilings}
\author{Carl Bracken \\
School of Mathematical Sciences\\
University College Dublin\\
Ireland}

\maketitle
\thispagestyle{empty}
\begin{abstract}
Site swap is a mathematical notation used by jugglers to communicate, create and study complex juggling patterns. Determining the number of possible site swap juggling patterns with respect to certain limiting parameters such as number of balls etc., is a problem that has been much studied and solved by many mathematicians. However, when the patterns have a throw height restriction (ceiling) the problem becomes difficult and is in general still open. In this article we derive some formulae for computing the number of possible juggling patterns with respect to certain ceiling types.
\end{abstract}

\section{Introduction}
Other than the fact that both mathematics and juggling have been with us for millenia, there seems to be little else connecting these two disciplines. Both have managed to develop in some form or another in almost all of the world's ancient cultures with hardly any interaction or overlap between these two groups. One exception would be Abu Sahl, who juggled glass bottles on the streets of $10^{th}$ century Baghdad before becoming a well known mathematician. It wasn't until the later part of the $20^{th}$ century that students of mathematics would have the opportunity to learn and practice the art of juggling. This opportunity came as juggling increased in popularity as a hobby and spread through the student societies of North American and later European universities. After some time it was noticed that many of those who attended the weekly juggling workshops where students of mathematics or physics. The correlation is not easily explained, but it seems that the type of people that enjoy juggling also enjoy mathematics. As enjoyment leads to practice and practice leads to excellence, there are many examples of mathematicians with exceptional juggling ability. Claude Shannon and Ron Graham are two well known juggling mathematicians and both have written papers on the mathematics of juggling \cite{CS}, \cite{RG}.

\section{Site Swap Juggling}
{\bf Preliminaries}\\
Site swap juggling notation is a concept that allows representation of idealised juggling patterns by a string of integers. The idea was developed independently by a number of people in the mid 1980s. We start by making a number of assumptions about the juggling patterns we wish to consider. 

\begin{itemize}
\item The types of objects being juggled are not specified and for convenience we call them balls.
The balls are fixed in number for any given pattern and each hand can hold, throw or catch only one at a time.

\item The juggling pattern is periodic with period $n$ so an action taken at time $t=a$ is taken again at $t=a+n$.

\item  Any number of hands can be used to juggle any pattern. If more than one hand is used then every hand throws with the same constant rhythm and in a strict ordering. Once a hand has made a throw then every other hand takes turns throwing at fixed length time intervals before this first hand throws again.

\item The paths of the hands and the balls are not considered, only the amount of time it takes for a thrown ball to be in a position that it can be thrown again.
\end{itemize}

We let $j^n_b$ denote a site swap juggling pattern with $b$ balls and period $n$. The pattern can be written as an $n$-tuple with a command $h_i$ at coordinate $i$ instructing the juggler how high to throw a ball at time $t=i$.
$$j^n_b=(h_0, h_1, h_2, . . . , h_{n-1})$$ 
Each $h_i$ specifies the amount of time it will take for the ball thrown at $t=i$ to travel through the air, be caught and be ready to be thrown again. The unit of time used is the interval between any two throws. This means that a ball thrown at time $t=i$ with height $h_i$ will be caught at time $i+h_i$ and can then be thrown again.
The assumption that every hand catches only one ball at a time means that, $$i+h_i \neq j+h_j \ mod \ n.$$
Less obviously, we have another restriction on the $h_i$, namely
$$\sum^{n-1}_{i=0}h_i =bn.$$
A full proof of the fact that $j^n_b$ is a site swap juggling pattern if and only if both of the above conditions are satisfied is given in \cite {RG}.\\

\noindent
{\bf Constructing a Site Swap Pattern.}\\
Once presented with an ordered set of $n$ integers it is easy to determine whether or not it represents a juggling pattern by checking whether it obeys both conditions above. If it is a legitimate pattern, we can also determine how many balls are required to juggle it by taking the average of the throw heights (a rearrangement of the second condition). However, it is also possible to directly construct an $n$-tuple that will obey both conditions by decomposing $j^n_b$ into three parts. That is, we let
$$j^n_b=P^n-Q^n+nB^n,$$
where $P^n=(p_0, p_1, ..., p_{n-1})$,  $Q^n=(0, 1, ..., n-1)$ and $B^n=(b_0, b_1, ..., b_{n-1})$. Each $p_i$ is a distinct element of the integers $\{0,1, ...,n-1\}$, each $b_i$ is a non negative integer and $\sum^{n-1}_{i=0}b_i=b$. It can be easily verified (as in \cite{RG}) that any vector constructed this way will be a site swap pattern and furthermore that every site swap pattern can be decomposed in this way.

We will now demonstrate the simplicity of this construction with an example. Let $n=b=4$, so we are constructing a 4 ball pattern with period 4. For $P^4$ choose $(1,2,3,0)$ and for $B^4$ choose $(1,1,1,1)$. Therefore,
$$j^4_4=(1,2,3,0)-(0,1,2,3)+4(1,1,1,1)=(5,5,5,1).$$
This is just one of the many new patterns discovered and now performed by jugglers since the introduction of site swap notation. Note that we can make any choice we wish for $P^n$, but $B^n$ must have $b_i \geq 1$ in any position where $P^n-Q^n$ is negative, otherwise the juggling pattern would contain a negative throw height. As height is a measurement of time, a negative height would mean that the ball was thrown backwards through time. This should be avoided. The requirement that $b_i$ be at least one in the positions where there are negative numbers in \ $P^n-Q^n$ \ is a crucial part of the counting argument in the next section. It is clear that $b_i$ need only be at least one in order to prevent negative numbers occurring in the pattern as every entry in $B^n$ is multiplied by $n$ while the lowest possible number in \ $P^n-Q^n$\  is \ $-(n-1)$.

Now that jugglers can construct all the possible site swap patterns for any period or number of balls they wish, the natural question is, when will this end? That is, how many site swap patterns are there?

\section{Counting Patterns}
 It's not too difficult to see that if we allow the period of the juggling patterns to extend to infinity then the number of possible site swap pattern will be infinite. So when we count the number of possible site swaps we put limitations on the parameters of the patterns. Let $J(n,b)$ denote the number of patterns of period $n$ and $b$ balls. In $\cite{RG}$ it is shown that $J(n,b)=(b+1)^n-b^n$. This result has been reproven in many different ways since its publication with significantly shorter proofs. In this section we will return to the method used in the original proof in order to derive a formula for the number of possible patterns limited by a maximum throw height as well as a specified number of balls and fixed period. 

The motivation for this is that when the above formula is applied to practical juggling, the count includes patterns that have unrealistically high throws. For example, if we count the number of 5 ball patterns with period 4 (both unremarkable numbers in this context), the formula gives 671 possible patterns. This count, however, includes patterns such as $(20,0,0,0)$. This coresponds to a ball spending approximately 4 seconds in flight, which would require a throw of about 20 meters. Even the world's best jugglers would find it difficult to throw a ball 5 meters in such a way that it can be comfortably caught and hence used in a pattern. A ceiling of 5 meters would roughly translate to a throw height of 10 or 11 for the $h_i$. Let $J(n,b,c)$ be the number of $j^n_b$ with each $h_i \leq c$. Ideally we would like to be able to count this for any choice of $c$. In this section we obtain a formula $J(n,b,c)$ whenever $c$ is of the form $an-1$ for any integer $a$. In the next section we will consider the number of patterns with $c \leq n-1$ and discuss its relation to the $rook(s,n)$ problem of Vardi \cite{IV}. First we recall some existing results from combinatorics.\\

\bigskip

\noindent
{\bf Eulerian numbers.}\\
Let $\mathcal{P}(n)$ denote the set of $n$-tuples $P^n$ as defined in section 2. Then the Eulerian number, denoted $E(n,k)$, is the number of all possible $P^n \in \mathcal P(n)$ such that $p_i < i$ for exactly $k$ values of $i$. This is not the original definition but it is shown in \cite{RG} to be equivalent. These numbers obey the recursive relation,
$$E(n,k)=(k+1)E(n-1,k)+(n-k)E(n-1,k-1).$$
Using this we can obtain the following array for small values of $n$. The number $E(n,k)$ is in the $k^{th}$ position on the $n^{th}$ row.

$$1$$
$$1 \ \ \ \ 1$$
$$1 \ \ \ \ 4 \ \ \ \ 1$$
$$1 \ \ \ \ 11\ \ \ \ 11 \ \ \ \ 1$$
$$1 \ \ \ \ 26 \ \ \ \ 66 \ \ \ \ 26 \ \ \ \ 1$$
$$1 \ \ \ \ 57 \ \ \ \ 302 \ \ \ \ 302 \ \ \ \ 57 \ \ \ \ 1$$

\noindent
We also have the identities
$$\sum^{n-1}_{k=0}E(n,k)=n!$$
and 
$$E(n,k)=E(n,n-k-1),$$
as well as the explicit representation
$$E(n,k)=\sum^n_{i=0}(-1)^i{n \choose i}[(k-i+1)^n-(k-i)^n].$$

\bigskip

\noindent
{\bf Number of ways to sum to an non negative integer.}\\
Let $B(n,b)$ denote the number of $n$-tuples of non negative integers with entries that sum to $b$. That is, the number of possible $B^n$ as defined in section 2. We have,
$$B(n,b)={n+b-1 \choose n-1}.$$
This result is well known and a simple proof can be found in \cite{SB}.
If we wish to determine the number of $B^n$ with each $b_i \leq a$ (denoted $B(n,b,a)$), then applying a standard inclusion-exclusion argument to the above identity yields,
$$B(n,b,a)=\sum^n_{i=0}(-1)^i{n \choose i}{n+b-1-i(a+1) \choose n-1}.$$

\bigskip

\noindent
{\bf Worpitzky's identity.}\\
The following equation, involving the Eulerian numbers, first appeared in \cite{JW} in 1881,
$$x^n=\sum^{n-1}_{k=0}E(n,k){x+k \choose n}.$$
This can be verified by applying an inductive argument on $x$ and using the recursive formula for $E(n,k)$.

\bigskip

\begin{thm}
The number of period $n$ site swap patterns with $b$ balls and ceiling of $an-1$ for any positive integer $a$ is given by 
$$J(n,b,an-1)=\sum^n_{i=0}(-1)^i{n \choose i}[(b-ia+1)^n-(b-ia)^n].$$
\end{thm}
Proof:\\
As every juggling pattern $j^n_b$ can be decomposed as $P^n-Q^n+nB^n$, we can compute the number of possible 
$j^n_b$ by taking the product of the number $P^n$ and the number of $B^n$. However, our choice of $B^n$ is restricted by the number of negative numbers that appear in $P^n-Q^n$, i.e., by our choice of $P^n$. If we could allow negative $h_i$ in $j^n_b$ and hence time travelling balls, then the number of patterns with each $h_i \leq an-1$ is $n!B(n,b,a-1)$.
This comes from $n!$ choices for $P^n$ times the $B(n,b,a-1)$ choices for the $B^n$ with each $b_i \leq a-1$ which only allows $j^n_b$ to be as high as $an-1$. If we wish to recount this without allowing negative $h_i$, then we have to have $b_i \geq 1$ in every position in $B^n$ where $P^n-Q^n$ is negative. If $P^n-Q^n$ has $k$ negative entries and we insist that in each of these $k$ positions $b_i \geq 1$, then the number of choices for $B^n$ will only be $B(n,b-k,a-1)$.
The number of $P^n-Q^n$ with $k$ negative entries is the same as the number of $p_i$ with $p_i < i$, i.e., it is the Eulerian number $E(n,k)$. Therefore summing over all Eulerian numbers for $k=0$ to $n-1$ and multiplying each one by the consequent number of choices for $B^n$ gives us all patterns without negative $h_i$.
Therefore we have
$$J(n,b,an-1)=\sum^{n-1}_{k=0}E(n,k)B(n,b-k,a-1).$$
This implies
$$J(n,b,an-1)=\sum^{n-1}_{k=0}E(n,k)\sum^n_{i=0}(-1)^i{n \choose i}{b+n-1-ia-k \choose n-1}$$
$$\ \ \ \ \ \ \ \ \ \ \ \ \ \ \ \ \ \ \ \ \ \ \ \ =\sum^n_{i=0}(-1)^i{n \choose i}\sum^{n-1}_{k=0}E(n,k){b+n-1-ia-k \choose n-1}.$$
Note that as $E(n,k)=E(n,n-k-1)$, we can replace $k$ with $n-k-1$ in the binomial coefficient while leaving it unchanged in $E(n,k)$ to obtain
$$J(n,b,an-1)=\sum^n_{i=0}(-1)^i{n \choose i}\sum^{n-1}_{k=0}E(n,k){b-ia+k \choose n-1}.$$
Now we apply the relation
$${x \choose n-1}={x+1 \choose n}-{x \choose n}$$
and obtain
$$J(n,b,an-1)=\sum^n_{i=0}(-1)^i{n \choose i}\sum^{n-1}_{k=0}E(n,k)[{b-ia+k+1 \choose n}-{b-ia+k \choose n}].$$
Finally we apply Worpitzky's identity to get
$$J(n,b,an-1)=\sum^n_{i=0}(-1)^i{n \choose i}[(b-ia+1)^n-(b-ia)^n]$$
and the proof is complete. \x

\section{Small ceilings and Vardi's rook problem}
From now on we will only consider ceilings that are less than the period. We shall refer to such ceilings as small, although they are only small relative to $n$ (which is only bounded by the jugglers memory and the audience's patience). As bounding height also bounds the maximum number of balls that can be juggled ($b \leq c$), we will also be considering $J(n,\ast ,c)$, where the `$\ast$' indicates that we are not fixing the number of balls in the count , i.e., $J(n,\ast,c)=\sum^c_{b=0}J(n,b,c)$.

From Theorem 1 we have $J(n,b,n-1)=\sum^n_{i=0}(-1)^i{n \choose i}[(b-i+1)^n-(b-i)^n]$. This is the explicit formula for the Eulerian number from section 2. This means 
$$J(n,b,n-1)=E(n,b)$$
and hence 
$$J(n,\ast ,n-1)=\sum^{n-1}_{b=0}J(n,b,n-1)=\sum^{n-1}_{b=0}E(n,k)=n!.$$
There is another way we could have arrived at these solutions. If we construct a juggling pattern with the decomposition method from section 2 and insist that each $h_i \leq n-1$, then the only places in $B^n$ we can have a nonzero entry are the positions where $P^n-Q^n$ are negative. This means there will always be only one choice for $B^n$ as we need $b_i=1$ whenever $p_i < i$ and the number of choices for $P^n$ will be $E(n,b)$. Therefore $E(n,b)$
will be the number of possible patterns with $b$ balls.
We find it interesting that by counting the number of possible juggling patterns in two different ways we can derive the explicit formula for $E(n,k)$.

\bigskip

Next we will consider ceilings that are smaller than $n-1$. When we count $J(n,\ast,n-1)$, there is only one choice for $B^n$ but we can use any $P^n$ as we have not specified the number of balls. If $c < n-1 $, there will again be only one choice for $B^n$ but the choices on $P^n$ will be limited also.
To count the number $J(n,\ast ,c)$ (whenever $c \leq n-1$) we need to find the number of $P^n$ such that $p_i-i \leq c \ mod \ n$ for all $p_i$ in $P^n$. An equivalent problem has already been studied by Vardi in \cite{IV}. It is a particular form of the rook placement with restrictions problem. In rook placement problems one has to place $n$ rooks on an $n \times n$ chessboard such that no rook can capture any other (i.e., one on every row and column). The task is to determine the number of possible arrangements when some other restrictions are introduced. When there are no extra restrictions the number of arrangements in $n!$.\\
We state (a version) of Vardi's problem: \\
Consider an $n \times n$ chessboard with the restriction that, for some fixed $s$ and any positive $t \leq s$ , a rook may not be put in column $t+i \ mod \ n $ when on row $i$, where the rows are numbered $0, 1,\ldots,n-1$. In \cite{IV} Vardi uses $rook(s,n)$ to denote the number of possible arrangements. He notes that $rook(1,n)$ is the number of derangements on $n$ symbols and $rook(2,n)$ is the solution to the married couples problem (see\cite{JT}). This is equivalent to counting the number of possible arrangements with no rooks on any $s$ specified and adjacent diagonals. 
By specifying the appropriate diagonals we can say that, this in turn is equivalent to the number of $P^n$ with each $p_i - i \leq n-s-1 \ mod \ n$. This implies that for $c \leq n-1$, we have
$$J(n,\ast,c)=rook(n-c-1,n).$$
Using the fact that both $rook(1,n)$ and  $rook(2,n)$ are known \cite{JT}, we can obtain two more formulae for $J(n,\ast ,c)$,
$$J(n,\ast,n-2)=\sum_{k=0}^n (-1)^k {n \choose k}(n-k)!$$
and
$$J(n,\ast,n-3)=\sum_{k=0}^n (-1)^k \frac{2n}{2n-k}{2n-k \choose k}(n-k)!.$$

\section{Closing Remarks}

In this article we have derived a number of expressions for the number of juggling patterns when a throw height limit is assumed. In particular, for juggling patterns of period $n$ we now have formulae for the number of patterns where the ceiling is $an-1$, $n-2$ or $n-3$. In the last two formulae we have not specified the number of balls and have instead summed over all possible patterns with all possible number of balls (including the $b=0$ pattern, which is neither difficult nor entertaining). This was done for the convenience of the mathematics and we would like to derive these formulae for a specified number of balls. In the $n-2$ case, this would be akin to partitioning the derangement numbers in the same way that the Eulerian numbers partition the factorials. It is not too hard to derive the three variable recursive relations for this problem, however an explicit solution seems difficult.

\end{document}